\documentclass[10pt]{amsart}
\newtheorem{thm}{Theorem}[section]
\newtheorem{cor}[thm]{Corollary}
\newtheorem{question}[thm]{Question}
\newtheorem{lem}[thm]{Lemma}
\newtheorem{prop}[thm]{Proposition}
\numberwithin{equation}{section}

\input epsf.tex

\begin{document}

\title{Ribbon tilings and multidimensional height functions}%
\author{Scott Sheffield}%
\address{Stanford University \\
Department of Mathematics Building 380 MC 2125 \\ Stanford, CA
94305
}%
\email{scott@math.stanford.edu}%

\thanks{This research was supported by a summer internship at Microsoft Research}%
\keywords{}%
\copyrightinfo{2001}
    {American Mathematical Society}

\begin{abstract}
We fix $n$ and say a square in the two-dimensional grid indexed by
$(x,y)$ has color $c$ if $x+y \equiv c \pmod{n}$.  A {\it ribbon
tile} of order $n$ is a connected polyomino containing exactly one
square of each color. We show that the set of order-$n$ ribbon
tilings of a simply connected region $R$ is in one-to-one
correspondence with a set of {\it height functions} from the
vertices of $R$ to $\mathbb Z^{n}$ satisfying certain difference
restrictions.  It is also in one-to-one correspondence with the
set of acyclic orientations of a certain partially oriented graph.

Using these facts, we describe a linear (in the area of $R$)
algorithm for determining whether $R$ can be tiled with ribbon
tiles of order $n$ and producing such a tiling when one exists. We
also resolve a conjecture of Pak by showing that any pair of
order-$n$ ribbon tilings of $R$ can be connected by a sequence of
local replacement moves.  Some of our results are generalizations
of known results for order-$2$ ribbon tilings (a.k.a. domino
tilings).  We also discuss applications of multidimensional height
functions to a broader class of polyomino tiling problems.

\end{abstract}
\maketitle
\section{Introduction}
\subsection{Ribbon tilings}
\begin{figure}
\begin{center}
\leavevmode \epsfbox[50 50 275 290]{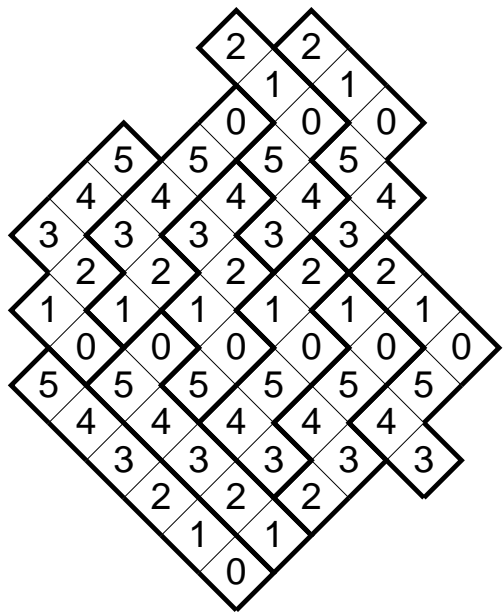}
\end{center}
\caption{A ribbon tiling of order $6$, rotated $45$ degrees
counterclockwise} \label{bigtiling}
\end{figure}

A square in the two-dimensional grid indexed by integers $(x,y)$
has color $c$ if $x+y \equiv c \pmod{n}$.  Two squares in the
$\mathbb Z^2$ lattice are {\it adjacent} if they share an edge. A
{\it ribbon tile} of order $n$ is a connected set of squares
containing exactly one square of each color.  (See Figure
\ref{bigtiling}.)  A ribbon tile can also be defined as a
connected sequence of $n$ squares, each of which comes directly
above or to the right of its predecessor.  A {\it region} $R$ is
any connected, finite subset of the squares of $\mathbb Z^2$. Two
order-$n$ ribbon tilings of $R$ are connected by a {\it local
replacement move} (or simply a {\it local move}) if one can be
obtained from the other by removing a single pair of ribbon tiles
and adding another pair in its place.

\begin{figure}
\begin{center}
\leavevmode \epsfbox[30 30 385 215]{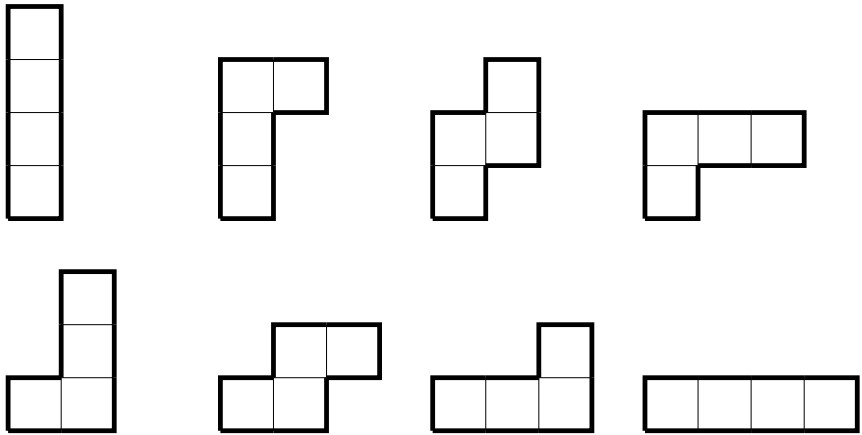}
\end{center}
\caption{Eight shapes for order-$4$ ribbon tiles}
\label{tetrominos}
\end{figure}

Labeled ribbon tilings of Young tableaux (also called {\it rim
hook tableaux}) have connections with symmetric function theory,
symmetric group representations, and monochromatic increasing
subsequences in colored permutations. (See \cite{B}, \cite{FS},
\cite{JK}, \cite{Ro}, and \cite{SW}.) Pak was the first to
consider ribbon tilings of more general regions \cite{P}.

A ribbon tile of order $n$ has one of $2^{n-1}$ possible shapes.
(See Figure \ref{tetrominos}.)  If $\epsilon$ represents one of
the $2^{n-1}$ shapes of a ribbon tile of size $n$ and $\alpha$ is
an order-$n$ ribbon tiling of a simply connected region, we denote
by $a_{\epsilon}(\alpha)$ the number of tiles in $\alpha$ of shape
$\epsilon$.

In \cite{P}, Pak coined the term {\it tile invariant} to describe
linear combinations of the $a_{\epsilon}(\alpha)$ whose values
were the same for all tilings $\alpha$ of a region $R$.  Pak
showed that if it were proved that any two order-$n$ ribbon
tilings of a region $R$ could be connected by a sequence of local
replacement moves, then certain tile invariants of $R$ could be
deduced as a trivial consequence \cite{P}.  Pak used a version of
the {\it Schensted algorithm for rim hook tableaux}---a
correspondence between rim hook tableau and $n$-tuples of smaller
tableaux (known by Nakayama and Robinson \cite{Ro} \cite{JK},
rediscovered by Stanton and White \cite{SW})---to prove that any
two order-$n$ ribbon tilings of a Young diagram could be
connected by a sequence of local moves, and that his tile
invariants must therefore hold for Young-tableau-shaped regions.
Pak used other techniques to extend the tile invariants to
row-convex regions and conjectured that both local move
connectedness and (consequently) the tile invariants could be
extended to arbitrary simply connected regions.

Without proving local move connectedness, Muchnik and Pak \cite
{MuP} and Moore and Pak \cite{MP} extended Pak's tile invariant
results to (respectively) order-$4$ and order-$n$ ribbon tilings
of general simply connected regions. Moore and Pak \cite{MP}
reiterated Pak's conjecture about local move connectedness for
general simply connected regions and cited private communication
from Thurston suggesting that the conjecture might be proved if
one could find a correspondence between ribbon tilings and
so-called {\it height functions} mapping the vertices of $R$ to
$\mathbb Z^n$.  Using Thurston's suggestion and a connection with
acyclic orientations, we prove the conjecture.

\begin{thm}
\label{connectedness} Any two order-$n$ ribbon tilings of a simply
connected region $R$ can be connected by a sequence of local
replacement moves.
\end{thm}

Here, we say $R$ is {\it simply connected} if $s \not \in R$
implies that there is no connected cycle of squares in $R$ that
encircles $s$ with a non-zero winding number; unless otherwise
specified, we will always assume $R$ is simply connected.

Our main ribbon tiling result relates the order-$n$ ribbon tilings
of a given region $R$ to the acyclic orientations extending a
certain partial orientation of a graph $G_R$. Pak's original paper
\cite{P} on ribbon tilings used a related construction (based on
earlier results in \cite{JK}, \cite{Ro}, and \cite{SW}) to prove
local replacement move connectedness for ribbon tilings of
Young-tableaux.  We will not use the results of these papers in
our presentation.

Before we can state the main theorem, we need to define several
terms. The first is a ``left of" relation, denoted $\prec$,
defined for both tiles and squares. In our pictures, we generally
rotate the set of tiles forty-five degrees so that $x+y$ increases
as we move upwards.  (See Figure \ref{bigtiling}.) Let $s_{x,y}$
be the square at position $(x,y)$.

We say $s_{x,y}$ is {\it directly left} of $s_{x',y'}$ if $x+y =
x'+y'$ and $x<x'$.  We say $s_{x,y}$ is {\it indirectly left} of
$s_{x',y'}$ if $|(x+y)-(x'+y')| = 1$ and $x \leq x'$ and $y \geq
y'$.  (Equivalently, $s_{x,y}$ is indirectly left of $s_{x',y'}$
if it is directly left of one of the four neighbors of
$s_{x',y'}$.)  If $s_{x,y}$ is directly or indirectly left of
$s_{x',y'}$, we write $s_{x,y} \prec s_{x',y'}$.

If $t$ is a tile and $s$ is a square, we write $s \prec t$ if $s
\prec s'$ for some $s' \in t$, and $t \prec s$ if $s' \prec s$ for
some $s' \in t$.  If $t_1$ and $t_2$ are two tiles in $\alpha$, we
write $t_1 \prec t_2$ if there exist $s_1 \in t_1$ and $s_2 \in
t_2$ with $s_1 \prec s_2$.  It is easy to verify that if $t_1$ and
$t_2$ are disjoint ribbon tiles, we cannot have both $t_1 \prec
t_2$ and $t_2 \prec t_1$.

We will also say a square $s_{x,y}$ is {\it higher} ({\it lower})
than $s_{x',y'}$ if $x+y>x'+y'$ ($x+y < x'+y'$). We say $x+y$ is
a square of {\it level} $c$ if $x+y = c$ and a tile $t$ has {\it
level} $c$ if $c$ is the level of its lowest square.  (A square's
color is equal to its level modulo $n$.)  An edge has {\it type}
$(i,i+1)$ if it separates squares of color $i$ and $i+1$, and a
tile has {\it type} $(i,i+1)$ if its lowest and highest squares
have color $i+1$ and $i$ respectively (when defining types, we
always add modulo $n$). Also, $a$ and $b$ are {\it comparable} by
$\prec$ if either $a \prec b$ or $b \prec a$.

Given a tiling $\alpha$, we define $t_{c,i}(\alpha)$ to be the
$i$th tile in $\alpha$, from left to right, of level $c$ (when
such a tile exists).  Clearly $t_{c,i}(\alpha) \prec
t_{c,i+1}(\alpha)$ whenever these two tiles are defined.

Let $B_R$ be the {\it boundary} of $R$, that is, the set of
squares that are not contained in $R$ but are adjacent to at least
one square in $R$.  We define the oriented graph $G_R(\alpha)$ to
be the graph whose vertices are the tiles in $\alpha$ and the
boundary squares of $R$.  Two such vertices $a$ and $b$ are
adjacent if they are comparable under $\prec$, and we say this
edge is {\it oriented} from $a$ to $b$ if $a \prec b$.

Some aspects of $G_R(\alpha)$ are independent of $\alpha$; using
height functions, we will prove the following result:

\begin{lem}
\label{sameGstructure} Any two tilings $\alpha$ and $\beta$ of $R$
necessarily contain the same number of tiles of each level. Thus,
$t_{c,i}(\alpha)$ is defined if and only if $t_{c,i}(\beta)$ is
defined.  Furthermore, if $c \equiv d \pmod{n}$, then
$t_{c,i}(\alpha) \prec t_{d,j}(\alpha)$ if and only if
$t_{c,i}(\beta) \prec t_{d,j}(\beta)$.  Also, if $b \in B_R$,
then $t_{c,i}(\alpha) \prec b$ if and only if $t_{c,i}(\beta)
\prec b$ .  Similarly, $b \prec t_{c,i}(\alpha)$ if and only if $b
\prec t_{c,i}(\beta)$.
\end{lem}

In other words, the map sending $t_{c,i}(\alpha)$ into
$t_{c,i}(\beta)$ is a canonical bijection between the tiles of
$\alpha$ and those of $\beta$, and this bijection preserves the
relation $\prec$ between tiles of the same type and between tiles
and boundary squares.

We define a tiling-independent partially oriented graph $G_R$ as
follows: all boundary squares in $B_R$ are vertices of $G_R$ and
we include a vertex called $t_{c,i}$ in $G_R$ if the tile
$t_{c,i}(\alpha)$ is defined for some tiling $\alpha$ (and hence
all such tilings) of $R$.  An edge is contained in $G_R$ if and
only if the corresponding edge is included in $G_R(\alpha)$ for
some tiling $\alpha$ (and hence all such tilings).  (Thus,
$t_{c,i}$ and $t_{d,j}$ are adjacent in $G_R$ if and only if
$|c-d| \leq n$; $t_{c,i}$ and $s_{x,y}$ are adjacent in $G_R$ if
and only if $c-1 \leq x+y \leq c+n$; and $s_{x,y}$ and $s_{x',y'}$
are adjacent in $G_R$ if and only $|(x+y) - (x'+y')| \leq 1$.)

Given a tiling $\alpha$, we can now think of $G_R(\alpha)$ as an
{\it orientation} on the graph $G_R$ (i.e., a way of assigning a
direction to each edge of $G_R$) in the obvious way: that is, an
edge $(t_{c,i},t_{d,j})$ is oriented from $t_{c,i}$ to $t_{d,j}$
in $G_R(\alpha)$ if and only if $t_{c,i}(\alpha) \prec
t_{d,j}(\alpha)$.  We define the orientation similarly for edges
involving two boundary squares or one tile and one boundary
square.

An edge involving two tiles of different types (i.e., some
$t_{c,i}$ and $t_{d,j}$ with $c \not \equiv d \pmod{n}$) is
called a {\it free edge} of $G_R$.  All other edges are called
{\it forced} edges of $G_R$.  By Lemma \ref{sameGstructure}, a
forced edge $(a,b)$ will be oriented in the same direction as
$G_R(\alpha)$ for every tiling $\alpha$.  Thus, we can think of
$G_R$ as being endowed with a partial orientation: its forced
edges are all oriented.  We say $G_R(\alpha)$ {\it extends} this
partial orientation by assigning an orientation to each free edge
of $G_R$.

Two orientations differ by an {\it edge reversal} if they agree on
all but exactly one edge.  The following is our main structure
theorem about the set of order-$n$ ribbon tilings of $R$.

\begin{thm}
\label{thecorrespondence} If $R$ admits at least one order-$n$
ribbon tiling, then there is a one-to-one correspondence between
the order-$n$ ribbon tilings of $R$ and the acyclic orientations
extending the partial orientation of $G_R$.  Furthermore, local
replacement moves in the space of tilings correspond to edge
reversals in the space of orientations. To be precise, we make
four assertions:
\begin{enumerate}
\item For every order-$n$ ribbon tiling $\alpha$ of $R$, $G_R(\alpha)$
is an acyclic orientation extending the partial orientation of
$G_R$.
\item If $\alpha$ and $\beta$ are order-$n$ ribbon tilings then $\alpha=\beta$ if
and only if $G_R(\alpha) = G_R(\beta)$.
\item Every acyclic orientation $A$ of $G_R$ that extends the partial orientation of $G_R$
is equal to $G_R(\alpha)$ for some tiling $\alpha$.
\item Two order-$n$ ribbon tilings $\alpha$ and $\beta$ differ by a local replacement
move if and only if $G_R(\alpha)$ and $G_R(\beta)$ differ by an
edge reversal.
\end{enumerate}
\end{thm}

Given Theorem \ref{thecorrespondence}, we will now deduce the
connectedness of ribbon tilings under local replacement moves by
citing a corresponding result about acyclic orientations.

Define the distance $d(A,B)$ between orientations $A$ and $B$ of
a graph $G$ to be the number of edges on which they differ. Then
it is well known that $A$ can be transformed into $B$ with
$d(A,B)$ edge reversals in such a way that each intermediate step
is also an acyclic orientation.  (See \cite{E} for a more general
result.)  Clearly, each intermediate step agrees with $A$ and $B$
on all edges on which $A$ and $B$ agree. It follows that if
$\alpha$ and $\beta$ are tilings, then the intermediate steps of
a length $d(G_R(\alpha), G_R(\beta))$ path of acyclic orientations
connecting $G_R(\alpha)$ and $G_R(\beta)$ must all be extensions
of the partial orientation of $G_R$. Thus, each of these
corresponds to an order-$n$ ribbon tiling of $R$.  Assuming
Theorem \ref{thecorrespondence}, we have proved Theorem
\ref{connectedness}.

See \cite{E}, \cite{St}, \cite{SZ} and the references therein for
more about the structure of the space of acyclic orientations. We
will also derive, as another consequence of Theorem
\ref{thecorrespondence}, an existence algorithm.

\begin{thm}
\label{thealgorithm} There is a linear-time (i.e., linear in the
number of squares of $R$) algorithm for determining whether there
exists an order-$n$ ribbon tiling of $R$ and producing such a
tiling when one exists.
\end{thm}

We begin in Section \ref{conwaysection} by reviewing the height
function theory of Conway-Lagarias and Thurston and showing that
it leads naturally to a construction of abelian height functions
for a variety of tiling problems.  In Section
\ref{ribbonsection}, we apply this theory specifically to the
ribbon tiling problem and prove our three main results: Lemma
\ref{sameGstructure}, Theorem \ref{thecorrespondence}, and
Theorem \ref{thealgorithm}. In Section \ref{generalsection} we
discuss generalizations of our abelian height function
constructions to tilings that do not appear to have the same
acyclic orientation characterization that ribbon tilings have.
Finally, in Section \ref{opensection} we present a number of open
problems in the theory of ribbon tilings and general abelian
height functions.

\section{Abelian Conway-Lagarias Thurston Height Functions}
\label{conwaysection} Our height function construction for ribbon
tilings makes use of a general technique developed by Conway and
Lagarias \cite{CL} and Thurston \cite{T} for analyzing tiling
problems.  A {\it polyomino} is a finite, connected subset of the
squares of $\mathbb Z^2$ whose complement is also connected.  For
simplicity, we will review the theory only for {\it polyomino
tilings}, although similar techniques apply if one replaces the
squares of the grid with the faces of any planar graph (e.g., the
hexagonal lattice).  (More general expositions---which invoke the
language of Cayley complex homology---can be found in \cite{Sch},
\cite{S}, and \cite{GP}.)

A {\it tile set} $T$ is a (usually infinite) set of finite,
simply-connected subsets of the squares of $\mathbb Z^2$.
Usually, we assume $T$ is translation invariant (i.e., $t \in T$
implies that $\{ s_{x,y}|s_{x+x_0, y+y_0}\in t \} \in T$ for all
$x_0, y_0 \in \mathbb Z$) but this need not be the case. The
first step in the Conway-Lagarias construction is to choose a map
$\phi$ from the oriented edges of the plane into a group $G$ in
such a way that $\phi(a) = -\phi(a^{-1})$ for all oriented edges
$a$; here $a^{-1}$ is the same edge as $a$ with opposite
orientation, and (since we will eventually restrict our attention
to abelian groups) we write the group operation additively. We
also require that the sum of $\phi(e)$ over the
clockwise-oriented edges $e$ on the boundary of any tile in $T$
(added in clockwise order) be equal to the identity in $G$. These
requirements are called {\it tile relations}.  Once the group is
chosen, given any tiling $\alpha$, we can construct a {\it height
function} $h_{\alpha}$ mapping at least some of the vertices of
$R$ into $G$ as follows.

First, choose a reference vertex $v_0$ on the boundary of $R$ and
require that $h_{\alpha}(v_0)$ be the identity.  Then if $v$ is
any other vertex of $R$ that is on the boundary of a tile,
$h_{\alpha}(v)$ is the ordered product of the oriented edges
along a path from $v_0$ to $v$ that does not cross any tiles. (If
the tiles are simply connected, such a path will always exist.)
Using the tile relations, it is not hard to prove that the value
$h_{\alpha}(v)$ is independent of the path chosen.

We can also always assume that the group $G$ is generated by the
values of $\phi(e)$ where $e$ is an edge in the $\mathbb Z^2$
lattice.  (Otherwise, replace $G$ by the subgroup generated by
these values.)  One can obtain the largest possible group $G$ that
is generated by $\phi(e)$ as follows:  Let $G$ be the group
generated by the oriented edges of $\mathbb Z^2$ modulo the
relations $\sum e = 0$ for clockwise sums of clockwise-oriented
boundary edges of a tile and $e+e^{-1}=0$ for all $e$.  Then let
$\phi(e)$ be the image of the edge $e$ in the group $G$.

The height function constructed in this way is in a sense
``maximally informative'' since any other height function is a
quotient of this one.  Unfortunately, this group is unwieldy in
practice and can be difficult to analyze. On the other hand, if
one merely seeks the maximally informative {\it abelian} group,
the analysis becomes much simpler.

Suppose $G$ is an additive abelian group and $\phi$ is a map from
the oriented edges of $\mathbb Z^2$ to $G$.  If $t$ is any
polyomino tile in $\mathbb Z^2$, we denote by $\phi(t)$ the sum of
the clockwise-oriented edges of $t$.  Observe that $$\phi(t) =
\sum_{s \in t} \phi(s),$$ where the sum is over squares $s$ that
are contained in the tile $t$.

Now, suppose that $p = (v_0, v_1, \ldots, v_k)$ is a
non-self-intersecting path of vertices in $R$ connecting the
reference vertex $v_0$ to a vertex $v_k$ on the boundary of some
tile in $\alpha$. We say that $p$ {\it crosses} a tile $t$ in
$\alpha$ if for some $0\leq i \leq k-1$, the edge $(v_i, v_{i+1})$
lies between two squares in $t$. If $s$ is a square in $R$ and
$0\leq i< j \leq k$, we say $p$ {\it passes $s$ on the left
between vertices $i$ and $j$} if:
\begin{enumerate}
\item The vertices $v_i$ and $v_j$ are on the boundary of a tile $t$
and all edges between $v_i$ and $v_j$ on $p$ are interior edges of
the tile $t$.
\item The square $s$ is contained in the polygon formed by the edges in $p$ between
$v_i$ and $v_j$ together with the counter-clockwise half of the
boundary of $t$ connecting $v_i$ to $v_j$.
\end{enumerate}
Note that a long path could pass a square $s$ on the left several
times if $s$ belongs to a tile that the path crosses several
times.

\begin{figure}
\begin{center}
\leavevmode \epsfbox[50 50 275 290]{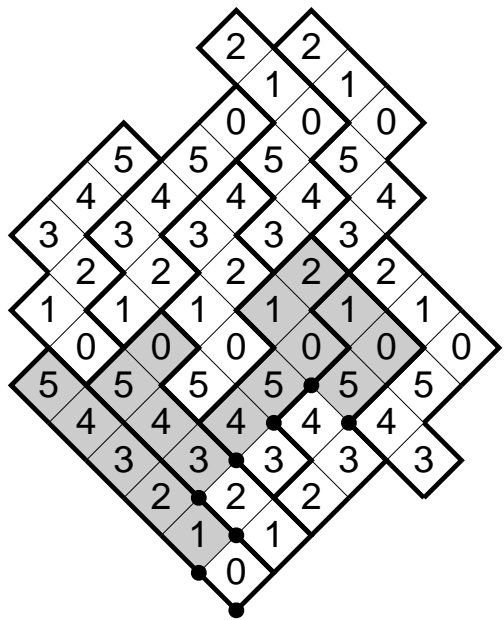}
\end{center}
\caption{This path of vertices (beginning with the lowest boundary
vertex) passes each shaded square on the left once.}
\label{shadedbigtiling}
\end{figure}

We define $\phi(p) = \sum_{i=0}^{k-1} \phi(v_i,v_{i+1})$ and
observe the following:

\begin{prop}
The height function value at $v_k$ can be written as
$$h_{\alpha}(v_k) = \phi(p) + \sum \phi(s),$$ where the latter sum
is over all squares (counted with multiplicity) that $p$ passes on
the left.
\end{prop}

To see this, let $p'$ be a path from $v_0$ to $v_k$ that is the
same as $p$ except that each interval that crosses a tile $t$ is
replaced by the path around the boundary of $t$ on the left side
of $p$.  The reader may easily check that $\phi(p') - \phi(p)$ is
equal to the sum of $\phi(s)$ over all squares (counted with
multiplicity) that $p$ passes on the left.  By definition,
$h_{\alpha} (v_k) = \phi(p')$, so the proposition follows.

Since $\phi(p)$ is independent of the particular tiling $\alpha$,
we have a simple interpretation of the height function
$h_{\alpha}(v)$: it tells us the sum of $\phi(s)$ over the squares
$s$ passed on the left by a path from $v_0$ to $v$.  For our
height function to be maximally informative, $G$ must be (or at
least contain) the abelian group generated by $\phi(s)$ subject
to only the tiling relations.  We refer to this group as the {\it
canonical abelian height space} (also known as the {\it color
group}) of a set of tiles $T$.

We say two squares $s_1$ and $s_2$ have the same {\it color} if it
can be deduced from the tiling relations that
$\phi(s_1)=\phi(s_2)$. (In other words, $s_1-s_2$ is zero in the
group generated by the squares of $\mathbb Z^2$ modulo the
relations $\sum_{s \in t} s = 0$ for all $t \in T$.)

This notion is best illustrated using our primary example: ribbon
tiling.  Let $s_{x,y}$ denote the square at position $(x,y)$ in
the lattice, and consider two overlapping horizontal ribbon tiles
of length $n$, one of which has leftmost square $s_{x,y}$ and the
other $s_{x+1,y}$.  The corresponding tile relations tell us that
$\sum_{i=0}^{n-1} \phi(s_{x+i,y}) = 0$ and
$\sum_{i=1}^{n}\phi(s_{x+i,y}) = 0$. And this implies that
$\phi(s_{x,y}) = \phi(s_{x+n,y})$.  Thus, $s_{x,y}$ and
$s_{x+n,y}$ have the same color.

Repeating similar arguments with other pairs of ribbon tiles that
overlap on all but one square, one can deduce that squares at
$(x_1,y_1)$ and $(x_2,y_2)$ necessarily have the same color
whenever $x_1+y_1$ and $x_2+y_2$ are equivalent modulo $n$. Denote
by $e_i$ the value of $\phi(s_{x,y})$ when $x+y \equiv i
\pmod{n}$. Then the relation $\sum_{i=0}^{n-1} e_i = 0$ is enough
to guarantee that $\phi(t)=0$ for any ribbon tile $t$ of length
$n$.  Thus, if $G$ is the abelian group generated by $e_0,
\ldots, e_{n-1}$ subject to this one relation, $G$ is maximally
informative.

At this point, we know the value of $\phi(s) \in G$ for every
square $s$.  However, we still have to choose a function $\phi$ on
the edges of the $\mathbb Z^2$ lattice in such a way that
$\phi(s)$ has the desired values.  This is essentially a linear
algebra problem, and it is not hard to see that it always has at
least one solution.  (We return to this issue in the next
section.)

This ``maximally informative'' abelian height function
construction is unique up to two obvious transformations.  First,
change of basis: one can replace $G$ with any group $G'$ and the
$e_i$ with any elements of $G'$ satisfying the necessary
relations.  Second, translation by a tiling-independent function:
one can add to $\phi(e)$ any function $\phi'$ satisfying
$\phi'(s)=0$ for all $s$.  This amounts to adding a function
(independently of $\alpha$) to each height function $h_{\alpha}$.

In particular applications, we try to choose the basis for $G$
and the values of $\phi$ on the edges in $\mathbb Z^2$ in a way
that makes the tiling space easy to understand.  In particular,
if $T$ is translation invariant, it is usually convenient for
$\phi$ to have some translational symmetries as well.

\section{Ribbon Tilings}
\label{ribbonsection}
\subsection{Height Functions for Ribbon Tilings}

Although we determined in the last section that $\phi(s)= e_i$
when $s$ is a square of color $i$, we have considerable freedom in
choosing the values of $\phi$ on the edges in $\mathbb Z^2$.  We
will use this freedom to force $\phi$ to have certain symmetries.

First, we would like $\phi$ to be invariant under color-preserving
translations of $\mathbb Z^2$; this amounts to the requirement
that $\phi(e)$ be the same whenever $e$ is a vertical (or
similarly, horizontal) edge of type $(i,i+1)$ oriented with color
$i+1$ on the left.  We would also like $\phi$ at such an edge to
have the same value---call it $e_{i,i+1}$---independently of
whether the edge is horizontal or vertical.

Let $e_{0,1}, e_{1,2}, \ldots, e_{n-1,0}$ be basis elements of a
height space $G = \mathbb Z^n$.  Summing around a square of color
$i$ tells us that $e_i$ must be equal to $2e_{i,i+1} -
2e_{i-1,i}$. Since this clearly implies $\sum_{i=0}^{n} e_i = 0$,
we need not add any additional relations to $G$.  (Although the
canonical abelian height function space is the subgroup of $G$
generated by the values $e_i = 2e_{i,i+1} - 2e_{i-1,1}$, it will
be notationally convenient to define $\phi$ and our height
functions in this larger group.) We can now describe our height
function $h_{\alpha}$ for a tiling $\alpha$ with the following
rules:

\begin{enumerate}
\item Let $h_{\alpha}(v_0)=0$ for all tilings $\alpha$ of $R$, where $v_0$ is a
reference vertex on the boundary of $R$.
\item When moving from one vertex to a neighbor along an edge
separating $i$ and $i+1$ color squares with color $i+1$ on the
left, the height changes by $e_{i,i+1}$ if no tile of $\alpha$ is
crossed.
\item When moving from one vertex to a neighbor along an edge separating
$i$ and $i+1$ color squares with color $i+1$ on the left, the
height changes by $2e_{j,j+1}-e_{i,i+1}$ if the edge crosses a
tile in $\alpha$ of type $(j,j+1)$.
\end{enumerate}

We write $h^{i,i+1}_{\alpha}$ to mean the $e_{i,i+1}$ component of
the height function $h_{\alpha}$.  Of course, there is some
redundancy in using an $n$ dimensional height space instead of the
canonical $n-1$ dimensional height space; as currently defined,
$\sum_{i=0}^{n-1} h^{i,i+1}_{\alpha}$ is a function that does not
depend on the tiling $\alpha$.  In particular, if $n=2$, the
projection sending $ae_{0,1}+be_{1,0}$ to the integer $a-b$
transforms our two-dimensional $h_{\alpha}$ into the standard
one-dimensional domino-tiling height function without any loss of
information.

\subsection{Interpretations of Ribbon Tiling Height Functions}

We can understand what the height function means by looking at
its behavior along diagonals.  We say $p$ is an {\it $(i,i+1)$
diagonal} if:

\begin{enumerate}
\item The path $p=v_0,\ldots, v_k$ is a connected left-to-right path
made solely of edges of type $(i,i+1)$ that are incident to
squares in $R$.
\item Both $v_0$ and $v_k$ are on the boundary of $R$ and the extending edges
$(v_{-1},v_0)$ and $(v_k,v_{k+1})$ of type $(i,i+1)$ are {\it not}
incident to squares in $R$.
\end{enumerate}

Note that as one moves from left to right along an $(i,i+1)$
diagonal, $h^{j,j+1}_{\alpha}$ increases by two if we cross a a
tile of type $(j,j+1)$ and zero otherwise. Thus, in the diagonal
described above, $h^{j,j+1}_{\alpha} (v_i)$ is equal to
$h^{j,j+1}_{\alpha}(v_0)$ plus twice the total number of times the
path $v_0, \ldots, v_i$ crosses a tile of type $(j,j+1)$.  In
particular, since the value of $h_{\alpha}$ on the boundary of $R$
is independent of $\alpha$, we have the following.

\begin{lem}
\label{diagonalnumberindependence} If $j \not = i$, then the total
number of tiles of type $(j,j+1)$ crossed by the portion of an
$(i,i+1)$ diagonal between two consecutive boundary vertices $v_1$
and $v_2$ is equal to $\frac{1}{2}(h_{\alpha}^{j,j+1}(v_2)
-h_{\alpha}^{j,j+1}(v_1))$.  In particular, since $v_1$ and $v_2$
are boundary vertices, this number is independent of $\alpha$.
\end{lem}

\begin{figure}
\begin{center}
\leavevmode \epsfbox[20 30 360 210]{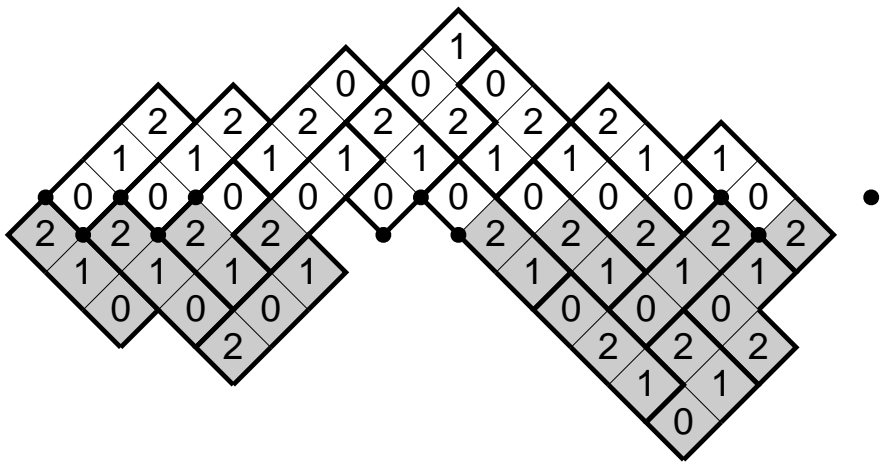}
\end{center}
\caption{Let $p$ be the diagonal of type $(2,0)$ that separates
the shaded squares from the unshaded squares. If $v_0$ is the
leftmost corner of the figure, then $h^{2,0}_{\alpha}(v_0)=0$. As
one moves from left to right along the diagonal,
$h^{2,0}_{\alpha}$ assumes, in order, the values
$0,1,2,3,4,5,4,5,4,5,6,7,8,7,8,7,8,7,8,9,10,9,10$.  The dots
indicate record vertices---one for the first occurrence of each
number between $1$ and $10$ and one for the first vertex outside
of $R$; the square immediately left of each record vertex is
either a boundary square or a square in a tile of type $(2,0)$.
The relation $\prec$ induces a total ordering on this set of four
boundary squares and seven type $(2,0)$ tiles. Every pair of
elements in this set corresponds to a forced edge of
$G_R(\alpha)$.} \label{diagonalexample}
\end{figure}

How does $h^{i,i+1}_{\alpha}$ change along an $(i,i+1)$ diagonal?
We know that moving from left to right along the diagonal, it
decreases by one when we cross a tile, increases by one otherwise.
Also, inspection shows that there can never be tiles crossing both
of two adjacent edges on a diagonal.  Thus, as one moves along the
diagonal from left to right, the value of $h^{i,i+1}_{\alpha}$
alternates between sequences of decreasing steps (of length at
most one) and sequences of increasing steps (of length one or
more). We say $v_a$ is a {\it record vertex} if $a \geq 1$, and
$h^{i,i+1}_{\alpha}(v_b)<h^{i,i+1}_{\alpha}(v_a)$ whenever $0 \leq
b < a$.  For notational convenience, we will also consider
$v_{k+1}$ (the vertex just right of the vertex $v_k$ along an edge
of type $(i,i+1)$) to be a record vertex.  (See Figure
\ref{diagonalexample}.)

Since each decreasing sequence in $h^{i,i+1}_{\alpha}$ along $p$
has length at most one, it is not hard to check that if $2 \leq a
\leq k$, then $v_a$ is a record vertex if and only if
$h^{i,i+1}_{\alpha}(v_a)>
h^{i,i+1}_{\alpha}(v_{a-1})>h^{i,i+1}_{\alpha}(v_{a-2})$.  Note
that the first edge $(v_0,v_1)$ and last edge $(v_{k-1}, v_k)$ of
a diagonal are always boundary edges of $R$, and hence represent
increasing steps.  In particular, $v_1$ is always a record vertex.

What does it mean for $v_a$ to be a record vertex?  First of all,
it means that the two adjacent edges $(v_{a-2},v_{a-1})$ and
$(v_{a-1},v_{a})$ are both {\it not} crossed by tiles; this in
turn implies that the square $s$ incident to both of them (which
we will call a {\it record square}) is either a square outside of
$R$ or a square in $R$ belonging to a tile of type $(i,i+1)$.
Since each of the integers $h^{i,i+1}_{\alpha}(v_0)+1,
h^{i,i+1}_{\alpha}(v_0)+2, \ldots, h^{i,i+1}_{\alpha}(v_k)$ must
occur as a record value exactly once, the number of record
vertices on the diagonal (not counting $v_{k+1}$) is equal to
$h^{i,i+1}_{\alpha}(v_k) - h^{i,i+1}_{\alpha}(v_0)$, and is hence
independent of $\alpha$.

Now, if $s$ is a record square outside of $R$ and $v_a$ is the
corresponding record vertex, the two diagonal edges preceding
$v_a$ must both be boundary edges of $R$ (unless one is an edge
outside of $R$); hence, $s$ will be a record square for every
tiling of $R$.  If $a \neq k+1$, then the value $h_{\alpha}(v_a)$
is also independent of the tiling.  Whenever $s$ is a record
square outside of $R$, the corresponding record vertex is called a
{\it boundary record vertex}; in particular, $v_1$ and $v_{k+1}$
are boundary record vertices.

Between any two consecutive boundary record vertices $v_a$ and
$v_b$, there are exactly $h^{i,i+1}_{\alpha}(v_b) -
h^{i,i+1}_{\alpha}(v_a)-1$ record vertices whose record squares
are squares of tiles of type $(i,i+1)$.  Whether the tile lies
above or below the diagonal clearly depends on the parity of $a$;
since $h^{i,i+1}$ alternates parity along the diagonal, this in
turn depends on the parity of the $h^{i,i+1}_{\alpha}$ at the
record vertex $v_a$.  We call $v$ a {\it lower record vertex} if
the corresponding record square lies below the diagonal and an
{\it upper record vertex} if it lies above the diagonal.  Since
for each $\alpha$, the record values $h^{i,i+1}_{\alpha}(v_0)+1,
h^{i,i+1}_{\alpha}(v_0)+2, \ldots, h^{i,i+1}_{\alpha}(v_k)$ must
be attained in the same order, we have the following lemma:

\begin{lem}
\label{diagonalordering} Suppose $q$ is the portion of a diagonal
$p$ of type $(i,i+1)$ between two consecutive boundary record
vertices $v_a$ and $v_b$. Then the number tiles of type $(i,i+1)$
with squares adjacent to $q$ that lie above or below the diagonal
is equal to $h^{i,i+1}_{\alpha}(v_b) - h^{i,i+1}_{\alpha}(v_a)-1$.
In particular, this value is the same for every tiling of $R$.
Furthermore, the order in which the entire diagonal $p$ encounters
boundary record vertices, lower record vertices, and upper record
vertices is also the same for every tiling of $R$ and can be
determined from the value of the height function on the boundary
vertices in $p$.
\end{lem}

Denote by $S_{\alpha}(p)$ the set of squares adjacent to $p$ that
are either members of tiles of type $(i,i+1)$ or squares in $B_R$.
(See Figure \ref{diagonalexample}.)  Note that if $s_1$ and $s_2$
are two record squares then $s_1 \prec s_2$ if and only if the
record vertex for $s_1$ occurs before that of $s_2$ in the
diagonal. Thus, we can restate the previous lemma:

\begin{lem}
\label{squarediagonalordering} Given any two tilings $\alpha$ and
$\alpha'$, $S_{\alpha}(p)$ and $S_{\alpha'}(p)$ contain the same
number of boundary squares, squares of color $i$ in type $(i,i+1)$
tiles, and squares of color $i+1$ in type $(i,i+1)$ tiles.  The
relation $\prec$ induces a total ordering on the squares in
$S_{\alpha}(p)$ as well as those in $S_{\alpha'}(p)$; if the
$m$th square in the ordering of $S_{\alpha}(p)$ is a boundary
square (or a square of color $i$ in a type $(i,i+1)$ tile, or a
square of color $i+1$ in a type $(i,i+1)$ tile), the same is true
of $S_{\alpha'}(p)$.
\end{lem}

Now we are ready to prove Lemma \ref{sameGstructure}, which we
restate here for convenience.

\begin{lem}
Any two tilings $\alpha$ and $\beta$ of $R$ necessarily contain
the same number of tiles of each level. Thus, $t_{c,i}(\alpha)$
is defined if and only if $t_{c,i}(\beta)$ is defined.
Furthermore, if $c \equiv d \pmod{n}$, then $t_{c,i}(\alpha)
\prec t_{d,j}(\alpha)$ if and only if $t_{c,i}(\beta) \prec
t_{d,j}(\beta)$.  Also, if $b \in B_R$, then $t_{c,i}(\alpha)
\prec b$ if and only if $t_{c,i}(\beta) \prec b$ .  Similarly, $b
\prec t_{c,i}(\alpha)$ if and only if $b \prec t_{c,i}(\beta)$.
\end{lem}

\begin{proof} Let $t_{c,1}(\alpha) \prec t_{c,2}(\alpha) \prec
\ldots \prec t_{c,k}(\alpha)$ be the tiles of type $(i,i+1)$ in
$\alpha$ whose lowest squares lie just above a given diagonal
$p$.  It is obvious from the definition of $t_{c,i}(\alpha)$ that
the ordering of these tiles with respect to each other is
independent of $\alpha$.  Similarly, let $t_{c-n,1}(\alpha) \prec
t_{c-n,2}(\alpha) \prec \ldots \prec t_{c-n, m}(\alpha)$ be the
tiles of type $(i,i+1)$ that lie just below the diagonal $p$.

Now, let $a_1 \prec a_2 \prec \ldots \prec a_r$ be the totally
ordered set of {\it all} the boundary squares and all of the $k+m$
tiles that are incident to $p$.  Since we know that the ordering
restricted to the set of tiles of form $t_{c,i}(\alpha)$ (or the
set of tiles of form $t_{c-n,i}(\alpha)$, or the set of boundary
squares) is independent of $\alpha$, it follows from Lemma
\ref{squarediagonalordering} that the total ordering $a_1 \prec
a_2 \prec \ldots \prec a_r$ is independent of $\alpha$.

Since any two tiles of the same type or any tile and any square
that are comparable by $\prec$ must be incident to {\it some}
common diagonal $p$, Lemma \ref{squarediagonalordering} covers all
possible cases in which the $\prec$ relation occurs between tiles
of the same type or between tiles and boundary squares.
\end{proof}

The reader may recall from the introduction that we used this
lemma used to define the partially oriented graph $G_R$ and to
justify our interpretation of $G_R{\alpha}$ as an extension of the
partial orientation of $G_R$.  This will be the topic of the next
section.

\subsection{The Acyclic Orientation Correspondence}

The purpose of this section is to prove Theorem
\ref{thecorrespondence}; we will state and prove each of its four
parts as a separate lemma.

\begin{lem} For every order-$n$ ribbon tiling $\alpha$ of $R$,
$G_R(\alpha)$ is an acyclic orientation extending the partial
orientation of $G_R$.
\end{lem}

\begin{proof} We must show that, given a tiling $\alpha$, there
exists no sequence $$a_0 \prec a_1 \prec a_2 \prec \ldots \prec
a_k \prec a_0,$$ where each $a_i$ is either a tile or a square in
$B_R$.

Since there can clearly be no cycle containing one or two
elements, we may assume $k \geq 2$.  Now, suppose a diagonal
incident to both $a_0$ and $a_k$ is also incident to some $a_i$
with $0 < i < k$. If $a_k \prec a_i$, then $a_i \prec a_{i+1}
\prec \ldots \prec a_k \prec a_i$ is a shorter cycle.  If $a_i
\prec a_k$, then $a_i \prec a_0$ (by transitivity, since $a_i$,
$a_k$, and $a_0$ are on the same diagonal and hence totally
ordered), and $a_0 \prec a_1 \prec \ldots \prec a_i \prec a_0$ is
a shorter cycle.

Thus, if $a_0 \prec a_1 \prec a_2 \prec \ldots \prec a_k \prec
a_0$ is a cycle of minimal length in $G_R(\alpha)$, then any
infinite diagonal that is incident to both $a_0$ and $a_k$ can be
incident to no other tile or square.   Denote by $p_z$ the
infinite diagonal between square levels $z$ and $z+1$ and let
$z'$ and $z''$ be the lowest and highest values of $z$ for which
$p_z$ is incident to both $a_0$ and $a_k$.  All tiles and squares
other than $a_0$ and $a_k$ must lie entirely below $p_{z'}$ or
entirely above $p_{z''}$. Since no tile or square below $p_{z'}$
can be comparable with any tile or square above $p_{z''}$, the
sequence $a_1 \prec a_2 \prec \ldots \prec a_{k-1}$ must lie
entirely above $p_{z''}$ or entirely below $p_{z'}$.  Suppose
without loss of generality that the former is the case and that
the lowest square in $a_0$ is at least as low as the lowest
square in $a_k$.  Then none of the elements $a_i$ below $p_{z'}$
is comparable with $a_k$, and this is a contradiction.
\end{proof}

\begin{lem} If $\alpha$ and $\beta$ are order-$n$ ribbon tilings
then $\alpha=\beta$ if and only if $G_R(\alpha) = G_R(\beta)$.
\end{lem}

\begin{proof} Obviously, the first statement implies the second.  Now,
suppose $G_R(\alpha) = G_R(\beta)$.  Since these orientations are
acyclic, there must exist at least one leftmost tile, i.e., there
exists a $t_{c,i}$ for which there exists no $t_{d,j}$ with
$t_{d,j}(\alpha) \prec t_{c,i}(\alpha)$ or $t_{d,j}(\beta) \prec
t_{c,i}(\beta)$. It follows that $t_{c,i}(\alpha)$ and
$t_{c,i}(\beta)$ must each contain the leftmost squares in $R$ in
the rows of level $c, c+1, \ldots, c+n-1$; thus, their positions
are determined and they are equal.  Similarly, we can now choose a
$t_{d,j}$ that is leftmost among the remaining set of tiles, and
its position is also completely determined.  Repeating this
argument, we conclude that $t_{c,i}(\alpha) = t_{c,i}(\beta)$ for
all $c$ and $i$ for which these values are defined.
\end{proof}

\begin{lem} \label{tilingfromorientation}
Every acyclic orientation $A$ of $G_R$ that extends the partial
orientation of $G_R$ can be written as $G_R(\alpha)$ for some
tiling $\alpha$.
\end{lem}

In the introduction, we defined $G_R$ by assuming that $R$
admitted at least one tiling $\alpha$ and restricting the
orientation on $G_R(\alpha)$ to the forced edges.  However, in our
later existence proofs, it will be necessary to have a definition
of $G_R$ and a version of Lemma \ref{tilingfromorientation} that
do not invoke {\it a priori} knowledge of the existence of
$\alpha$. To do this, we will need some definitions.

We say a tile vertex $t_{c,i}$ {\it crosses level $d$} if $c \leq
d \leq c+n-1$.  (Note that if $\alpha$ is a tiling, then $t_{c,i}$
crosses level $d$ if and only if $t_{c,i}(\alpha)$ contains a
square of level $d$.)  We define a partially oriented graph $G$ to
be a {\it ribbon tile graph consistent with $R$} if the following
hold:

\begin{enumerate}
\item The vertices of $G$ come from the set of tile vertices $t_{c,i}$
(with $c,i \in \mathbb Z$ and $i > 0$) and square vertices
$s_{x,y}$ (with $x,y \in \mathbb Z$).
\item Two tile vertices $t_{c,i}$ and $t_{d,j}$ are adjacent if and only if $|c-d| \leq n$.
A tile $t_{c,i}$ and square $s_{x,y}$ are adjacent if and only if
$c-1 \leq x+y \leq c+n$. Two squares $s_{x,y}$ and $s_{x',y'}$ are
adjacent if and only if $|(x+y) - (x'+y')| \leq 1$.
\item The graph $G$ is endowed with an orientation on its
{\it forced} edges (i.e. all edges involving a square or
involving two tile vertices $t_{c,i}$ and $t_{d,j}$ with $c
\equiv d \pmod{n}$).  No other edges of $G$ are oriented.
\item If $s_1$ and $s_2$ are square vertices, then the edge $(s_1,s_2)$
is oriented toward $s_2$ if and only if $s_1 \prec s_2$. Also,
whenever $i>0$ and $t_{c,i+1}$ is a tile vertex, $t_{c,i}$ is
also a tile vertex.  An edge of the form $(t_{c,i}, t_{c,j})$ is
oriented toward $t_{c,j}$ whenever $i<j$.
\item The set of square vertices of $G$ is precisely the set of boundary squares
of $R$.
\item The total number of tiles crossing a given level $c$ is precisely the
number of squares in $R$ of level $c$.
\item For every integer $c$ and every pair of boundary squares $b_1$ and $b_2$
of levels $c-1$, $c$, or $c+1$ satisfying $b_1 \prec b_2$, the
number of tile vertices $t_{d,j}$ in $G_R$ that cross level $c$
and satisfy $b_1 \prec t_{d,j} \prec b_2$ is precisely the number
of squares $s$ in $R$ of level $c$ satisfying $b_1 \prec s \prec
b_2$.
\end{enumerate}

Only the last three conditions depend on $R$; we will refer to
these as the {\it consistency} conditions.  Although it involves
overloading the symbol $\prec$, we will write $a \prec b$ whenever
the edge $(a,b)$ is a forced edge of $G$ oriented toward $b$.  If
$A$ is an orientation extending the partial orientation on the
forced edges of $G$, we write $a \prec_A b$ if the edge $(a,b)$ in
$G$ is oriented from $a$ to $b$ in $A$.

Note that up to this point, we have defined $\prec$ only on
actual tiles, not on the tile vertices $t_{c,i}$, which are
abstract place holders for tiles.  However, our expanded use of
$\prec$ is consistent in the following sense: if $G_R(\alpha)$ is
the ribbon tiling graph of a tiling $\alpha$ and the edge
$(t_{c,i}, t_{d,j})$ is forced, then $t_{c,i}(\alpha) \prec
t_{d,j}(\beta)$ if and only if $t_{c,i} \prec t_{d,j}$ in $G_R$.
Similarly, if $A=G_R(\alpha)$ and the edge is $(t_{c,i},t_{d,j})$
is not necessarily forced, then $t_{c,i}(\alpha)\prec
t_{d,j}(\alpha)$ if and only if $t_{c,i} \prec_A t_{d,j}$.
Similar results apply to edges involving squares.

If $R$ admits a tiling $\alpha$ then the graph $G_R$ as defined in
the introduction clearly satisfies the above list of requirements
and is consistent with $R$.  (In fact, $G_R$ is the only ribbon
tile graph consistent with $R$ in this case.)  Thus, the following
is a slightly more general version of Lemma
\ref{tilingfromorientation}:

\begin{lem} \label{tilingfromconsistentgraph}
Suppose $G_R$ is any ribbon tile graph that is consistent with
$R$. Then every acyclic orientation $A$ of $G_R$ that extends the
partial orientation on the forced edges of $G_R$ is equal to
$G_R(\alpha)$ for some tiling $\alpha$.  (In particular, if $R$
does not admit an order-$n$ ribbon tiling, then no such $G_R$ and
$A$ exist.)
\end{lem}

\begin{proof}
Suppose $t = t_{c,1}$ is a leftmost tile vertex in $G_R$. As in
the previous proof, we will choose the leftmost squares of each of
the levels $c, c+1, \ldots, c+n-1$ and declare them to be our
first tile $t_{c,1}(\alpha)$.  That there exists at least one
square in $R$ on each of these levels follows from the consistency
conditions of $G_R$.  However, we must verify that this set of
squares is actually a ribbon tile.  Since it contains one square
on each of levels $c, c+1, \ldots, c+n-1$, it is sufficient to
show that it is connected.

Suppose otherwise; then we can assume without loss of generality
that there exists some $i$ such that the leftmost square $s_i$ in
the row of level $c+i$ is left of the square $s_{i+1}$ in the row
of level $c+i+1$ but that these two squares are not adjacent.
Now, let $s_r$ ($s_l$) be the square of level $c+i+1$ that is
right of (left of) and adjacent to $s_i$.  Both $s_r$ and $s_l$
must be squares in $B_R$, since they are adjacent to $s_i$ and
they are left of $s_{i+1}$, the leftmost square of their level in
$R$.

Now, by the consistency conditions, there is exactly one tile $t'$
containing a square on level $c+i$ such that $s_l \prec t' \prec
s_r$; every other tile with a square on level $c+i$ must be right
of $s_r$. If $t \not = t'$, then $t' \prec s_r \prec t \prec t'$,
violating acyclicity. Hence $t = t'$. However, since $s_r$ is left
of every square of level $c+i+1$, one may deduce from the
consistency conditions that every tile vertex $t$ of $G_R$ that
crosses level $c+i+1$ must satisfy $s_r \prec t$. Thus $t \prec
s_r \prec t$, a contradiction, and it follows that the squares
chosen do in fact form a tile.

Now we form a new region $R'$ by removing the squares in this tile
from $R$.  We form a new graph $G_{R'}$ by removing $t_{c,1}$
from the vertex set of $G_R$---along with all boundary squares of
$R$ that are not boundary squares of $R'$---and adding all the
boundary squares of $R'$ that were not boundary squares of $R$.
(These are necessarily members of our removed tile.)  We define
the edges of $G_{R'}$ and their orientations in a new orientation
$A'$ in the most natural possible way.  Two vertices of $G_{R'}$
are adjacent whenever the corresponding edges in $G_R$ were
adjacent, and $A'$ is the orientation induced by $A$ on such
edges. Furthermore, if $s$ is one of the new boundary squares and
is in a row $c$, we declare it to be adjacent in $G_{R'}$ to
every tile $t$ that crosses rows $c-1$, $c$, or $c+1$ (with $s
\prec t$) and every boundary square $s_0$ of those levels (with
$s \prec s_0$ if and only if $s$ is directly or indirectly left
of $s_0$ and $s_0 \prec s$ when the opposite is true). Finally, if
$t_{c,j}$ was a vertex of $G_{R}$ with $j>1$, we relabel it
$t_{c,j-1}$ in $G_{R'}$.

It follows easily from the definition that the $G_{R'}$ thus
defined---together with the partial orientation on forced edges
induced by $A'$---is a ribbon tile graph consistent with $R'$.  If
we can further prove that $A'$ is necessarily acyclic, then Lemma
\ref{tilingfromconsistentgraph} will follow by induction.

Suppose $A'$ has a cycle; then this cycle must include a square
$s$ of $B_{R'}$ that was a member of the tile $t$---otherwise the
cycle would also be a cycle in $A$.  Also, although the cycle may
contain multiple squares from the removed tile, it must contain at
least one element that is not a square of $t$, since the squares
do not themselves contain any cycles.  Suppose the cycle contains
a sequence $s_1 \prec \ldots \prec s_k$ of former squares of $t$
but that the elements $a$ immediately preceding this sequence and
$b$ immediately following it are not squares of $t$.  Then we must
have had $a \prec t \prec b$ in $A$, and $A$ too must have had a
cycle.
\end{proof}

\begin{lem}
\label{localreplacementmove} Two order-$n$ ribbon tilings
$\alpha$ and $\beta$ differ by a local replacement move if and
only if $G_R(\alpha)$ and $G_R(\beta)$ differ by an edge reversal.
\end{lem}

\begin{figure}
\begin{center}
\leavevmode \epsfbox[18 40 390 220]{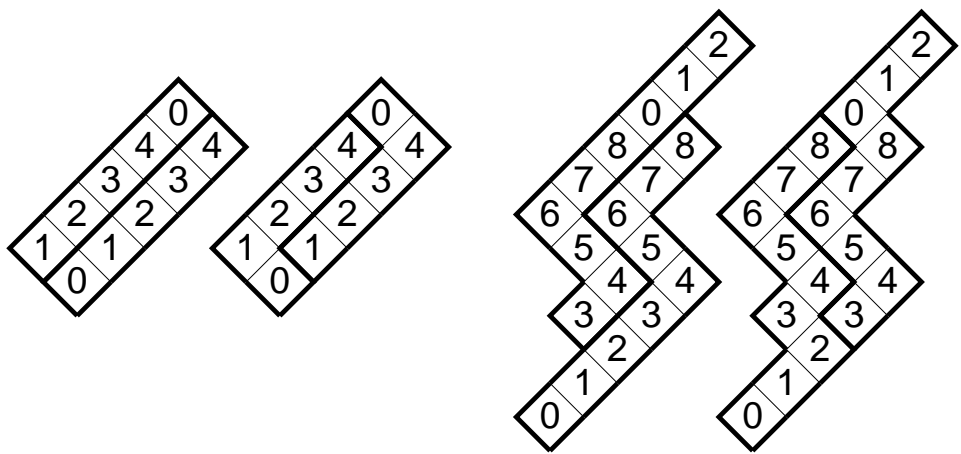}
\end{center}
\caption{Two examples of local replacement moves}
\label{localmoves}
\end{figure}

\begin{proof} First of all, we can determine whether a local replacement move
involving tiles $t_{c,i}(\alpha)$ and $t_{d,j}(\alpha)$ can occur
by applying the theory we have already developed to $R'$, the
region tiled by just these two tiles.  Now, if $R'$ admits two
distinct tilings, then $G_{R'}$ must admit two distinct
orientations that extend the orientation on the forced edges.  In
particular, $G_{R'}$ must contain at least one free edge, and
that can only be the edge $(t_{c,i}, t_{d,j})$.  For this edge to
be free, we must have $0<|c-d|<n$.

Without loss of generality, suppose $c < d$. Notice that if $d
\leq e \leq c+n-1$, then $R'$ contains two squares of level $e$.
Inspection shows that if tile $t_{c,i}(\alpha)$ contains the
leftmost square of one of these same-level pairs, it must contain
the leftmost square of every such pair.  (See Figure
\ref{localmoves}.)  Moreover, if a tile contains one square on a
row higher than $c+n-1$ (lower than $d$), it must contain every
square with this property. Thus, a local replacement move must
consist of swapping the squares in each same-level pair. Such a
swap produces a new pair of ribbon tiles if and only if the
single square on row $d-1$ (resp., $c+n$) is adjacent to both of
the squares in row $d$ (resp., $c+n-1$). Inspection shows that
whenever there is no square $s$ with either $t_{c,i}(\alpha)
\prec s \prec t_{d,j}(\alpha)$ or $t_{d,j}(\alpha) \prec s \prec
t_{c,i}(\alpha)$, then a local move can in fact occur.

Now, if $\beta$ is the tiling obtained by applying such a move to
$\alpha$, and $t_{c,i}(\alpha) \prec t_{d,j}(\alpha)$ then
inspection shows that we must have $t_{d,j}(\beta) \prec
t_{c,i}(\beta)$; it is also easy check that if $a$ is any tile or
square outside of $R'$, then $a \prec t_{c,i}(\alpha)$
 if and only if $a \prec t_{c,i}(\beta)$.  Similarly, $t_{c,i}(\alpha) \prec a$ if and
only if $t_{c,i}(\beta)\prec a$.  It follows that if $\alpha$ and
$\beta$ differ by a local replacement move, $G_R(\alpha)$ and
$G_R(\beta)$ do in fact differ by the reversal of exactly one
edge: the edge $(t_{c,i},t_{d,j})$.

For the converse, note that if the orientation of an edge
$(t_{c,i}, t_{d,j})$ of the orientation $G_R(\alpha)$ can be
reversed without creating a cycle, we know there must be no
vertex $a$ of $G_R(\alpha)$ with either $t_{c,i}(\alpha) \prec a
\prec t_{d,j}(\alpha)$ or $t_{d,j}(\alpha) \prec a \prec
t_{c,i}(\alpha)$; we have already shown that this is a sufficient
condition for there to exist a local move involving the tiles
$t_{c,i}(\alpha)$ and $t_{d,j}(\alpha)$, and so $G_R(\beta)$ must
be the orientation obtained by reversing the edge.
\end{proof}

\subsection{The Existence Algorithm}
The main result of this section will be a linear-time algorithm
for determining whether a simply connected region $R$ has an
order-$n$ ribbon tiling and constructing such a tiling when one
exists.  To determine existence, we will make use of the
well-known fact (which we later prove) that there exists an
acyclic orientation $A$ of a graph $G$ that extends a given
partial orientation of $G$ if and only if the partial orientation
is itself acyclic. This implies the following corollary of Lemma
\ref{tilingfromconsistentgraph}:

\begin{cor}
Given $R$, there exists an acyclic partially oriented ribbon
tiling graph $G_R$ that is consistent with $R$ if and only if $R$
admits an order-$n$ ribbon tiling.
\end{cor}

Using this fact, we divide our algorithm into three steps.

\begin{enumerate}
\item We either produce an acyclic ribbon tile graph $G_R$
that is consistent with $R$ or prove that no such graph exists
(and hence $R$ admits no order-$n$ ribbon tiling).
\item If we succeed in producing $G_R$ in the previous step, we
describe a total ordering on the vertices of $G_R$ that induces
the partial orientation of $G_R$ on the forced edges of $G_R$.
\item Given this totally ordering on the vertices of $G_R$, we
describe a tiling of $R$.
\end{enumerate}

Now, we know that if there is an acyclic ribbon tile graph $G_R$
that is consistent with $R$, we can use it to generate
$G_R(\alpha)$ using the construction of Lemma
\ref{tilingfromconsistentgraph}, and the partial orientation on
$G_R$ will thus be the restriction of $G_R(\alpha)$ to the forced
edges in $G_R$.  Thus, in the construction that follows, we can
make use of lemmas proved earlier in the paper under the
assumption that $R$ had a tiling.  If $R$ does have a tiling, our
assumption will have been correct.  If $R$ does not have a
tiling, these assumptions will lead to contradictions, and these
will enable us to conclude that no tiling exists.

Now we begin to construct $G_R$.  Let $S(c)$ be the number of
squares of level $c$ and let $T(c)$ be the number of tiles whose
lowest square is $c$.  Clearly $T(c) = S(c) - \sum_{i=1}^{n-1}
T(c-i)$, and is thus easily computed for each $c$ in linear time.
(If it is not the case that $T(c) \geq 0$ whenever $S(c)\not = 0$
and $T(c) = 0$ whenever $S(c) = 0$, we can conclude that no tiling
exists.)  For each $c$ for which $T(c) \not = 0$, we include
$t_{c,i}$ for $1 \leq i \leq T(c)$ as tile vertices of $G_R$. We
also include all the boundary squares of $R$ as square vertices
of $G_R$.  This $G_R$ has an edge whenever the corresponding
tiles/squares are on levels that would make them comparable by
$\prec$, and we require that $t_{c,i} \prec t_{c,i+1}$ whenever
both of these vertices are in $G_R$.

We would now like to determine an orientation on the forced edges
of $G_R$ that will make $G_R$ a consistent ribbon tiling graph.
We will not store in memory all of the edges of $G_R$ (as the
number of such edges could be quadratic in the size of $R$);
instead, we will observe that the partial orientation on $G_R$
induces a total ordering $a_0 \prec a_1 \prec \ldots \prec a_k$
on the set of squares and tiles of a given type that are all
incident to a common diagonal.  Since every forced edge can be
derived from such an ordering on some diagonal, it will thus be
sufficient for us to store in memory the actual orderings
corresponding to each such type/diagonal pair, and this requires
only a linear amount of memory.

Given an infinite diagonal $p$ between square rows of levels $c-1$
and $c$ and some $j \in \{1,2,\ldots,n-1\}$, let $\mathcal
S_{j,p}$ be the set of tiles in $G_R$ of the form $t_{c-j,i}$
together with squares in $G_R$ on levels $c-1$ and $c$.  Since
the elements in $S_{j,p}$ share a common diagonal and all tiles
in this set are of the same type, each pair of elements in
$S_{j,p}$ is connected by a forced edge and hence $G_R$ induces a
total ordering on this set.  Now, we can compute the height
function $h$ on the boundary of $R$, and using $h$, we can deduce
this ordering explicitly from Lemma
\ref{diagonalnumberindependence} as follows.  Let $k = c-j-1
\pmod{n}$. Given consecutive boundary vertices $v_1$ and $v_2$
along the diagonal, we know that the number of tiles of the form
$t_{c-j,i}$ (and hence of type $(k,k+1)$) that cross $p$ between
consecutive boundary vertices $v_1$ and $v_2$ must be
$\frac{1}{2}(h^{k,k+1}(v_2) - h^{k,k+1}(v_1))$.  Clearly, if
$b_1$ is a square incident to the diagonal $p$ immediately before
$v_1$ and $b_2$ is a square incident to $p$ immediately after
$v_2$, the number of tiles $t$ of the form $t_{c-j,i}$ satisfying
$b_1 \prec t \prec b_2$ must be exactly
$\frac{1}{2}(h^{k,k+1}(v_2)-h^{k,k+1}(v_1))$.  Thus, we can
compute the number of tiles of $S_{j,p}$ in between each pair of
boundary squares, and this completely determines the ordering.

Similarly, given an infinite diagonal $p$ between square rows
$c-1$ and $c$, let $\mathcal S_{p,0}$ be the set of tiles in $G_R$
of the form $t_{c,i}$ or $t_{c-n,i}$ together with the squares in
$G_R$ on levels $c-1$ and $c$.  Each pair of elements in
$\mathcal S_{p,0}$ is connected by a forced edge, and thus the
partial orientation on $G_R$ must induce a total ordering on this
set.  By Lemma \ref{diagonalordering} and Lemma
\ref{squarediagonalordering}, we can deduce this total ordering
from the value of the height function on the boundary vertices of
$R$ in $p$.

Since every forced edge involves two elements along the same
diagonal, we can compute the orientation of each forced edge by
repeating this process for each diagonal.  Now we check that if
$G_R$ is acyclic, it is necessarily a consistent ribbon tiling
graph.  It is enough to check the consistency conditions.  By
construction, $G_R$ contains all boundary squares, and by the
definition of $T(c)$, the total number of tiles crossing each
level $c$ is the number of squares in $R$ of level $c$.  For the
final consistency condition, suppose $c$ is given and $b_1$ and
$b_2$ are boundary squares of types $c-1$, $c$, or $c+1$ with $b_1
\prec b_2$.

Then let $p_1$ be the diagonal above the row of squares of level
$c$ and let $p_2$ the diagonal below these squares.  One easily
checks that both $p_1$ and $p_2$ contain edges in each of $b_1$
and $b_2$. It suffices to check the condition under the assumption
that there exists no boundary square $b_3$ incident to $p$ with
$b_1 \prec b_3 \prec b_2$, since otherwise we can check the
condition separately for the pair $(b_1,b_3)$ and the pair
$(b_3,b_2)$. Say a vertex has level $c$ if it is a corner of
squares of type $c-1$, $c$, and $c+1$, and let $q$ be the portion
of the diagonal connecting the last level $c$ vertex of $b_1$ to
the first level $c$ vertex of $b_2$.

Then it follows by applying Lemma \ref{diagonalnumberindependence}
to the diagonal $p_1$ (or to $p_2$ if $p_1$ has type $(j,j+1)$)
that the number of tile vertices $t$ of type $(j,j+1)$ crossing
$c$ and satisfying $b_1 \prec t \prec b_2$ must be equal to
$\frac{1}{2}(h^{j,j+1}(v_2) - h^{j,j+1}(v_1))$.  Write
$\overline{h} = \sum_{j=0}^{n-1} h^{j,j+1}$, and note that the
total number of such tile vertices must then be
$\frac{1}{2}(\overline{h}(v_2) - \overline{h}(v_1))$.  However,
this function is tiling-independent and in fact {\it always}
increases by one as we move between squares of color $j$ and $j+1$
with $j+1$ on the left.  It follows that if the vertices of the
grid are labeled with coordinates $(u,v)$, then
$\overline{h}(u,v)$ will be equal to $u-v$ plus a constant.  In
particular, this implies that $\frac{1}{2}(\overline{h}(v_2) -
\overline{h}(v_1))$ is equal to the number of squares on level $c$
directly between vertices $v_1$ and $v_2$, and this is easily seen
to be the number of level $c$ squares $s$ in $R$ satisfying $b_1
\prec s \prec b_2$.

\begin{figure}
\begin{center}
\leavevmode \epsfbox[45 95 395 220]{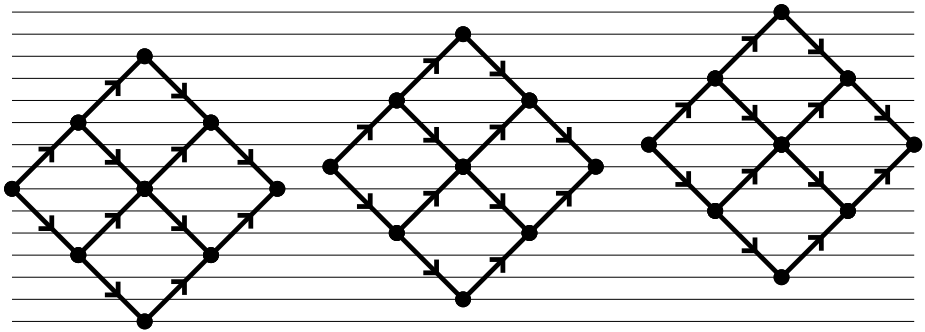}
\end{center}
\caption{A graphical representation of $G_R$ when $n=3$ and $R$
is the $9 \times 9$ grid.  Each of the dots represents one of the
$27$ tile vertices of $G_R$ (square vertices are not shown), and
the $15$ horizontal lines correspond to the $15$ levels on which
tiles occur.  All essential forced edges between tiles of the
same type are shown with orientations.  Two vertices in this
figure are adjacent in $G_R$ if and only if the difference
between their levels is at most three.} \label{interlockinggrids}
\end{figure}

The only remaining step is to determine whether our partial
orientation of $G_R$ is indeed acyclic.  There are many of ways to
do this in time linear in the number of edges of $G_R$.  We say a
forced edge $(a,b)$ of $G_R$ is {\it essential} if $a \prec b$
and there exist no forced edges $(a,c)$ and $(c,b)$ with $a \prec
c$ and $c \prec b$---one easily checks that the number of
essential forced edges of $G_R$ is linear in the size of $R$ and
that any acyclic orientation of $G_R$ that extends the partial
orientation on the essential forced edges must extend the partial
orientation on all forced edges.  (See Figure
\ref{interlockinggrids}.)  The following algorithm for determining
acyclicity (know as the {\it topological sort}) is well-known,
and we will apply it to the graph formed by the vertices of $G_R$
and the essential forced edges of $G_R$.

Choose a starting vertex $a_0$ from the set of vertices of $G_R$
and attempt to construct an increasing path $a_0 \prec a_1 \prec
a_2 \ldots$ until we either reach a point we have already visited
(and we declare that there exists a cycle) or we reach a vertex
$a_k$ that is not left of any of its neighbors. Clearly, $a_k$
cannot be part of a cycle, so we remove it and begin searching
again with $a_{k-1}$ (if $k=0$, we choose a new $a_0$ arbitrarily
from among the remaining vertices of $V$).  We continue this
process until we have removed all the vertices of $V$ or we have
found a cycle.

In the former case, if we list the elements of $V$ in the order
in which they were removed from $V$, we have some ordering $a_0,
a_1, a_2, \ldots, a_k$ in which $i < j$ implies $a_i \not \prec
a_j$ (since $a_i$ could otherwise not have been removed before
$a_j$).  Thus, we can define an acyclic orientation $A$ on $G_R$
by simply saying that whenever $(a_i, a_j)$ is an edge in $G_R$,
$a_j \prec_A a_i$.  Let $\alpha$ be the tiling for which
$G_R(\alpha) = A$. Then if $t_{c,i}$ is the last tile that occurs
in our ordering, we know that $t_{c,i}(\alpha)$ contains the
leftmost squares of rows $c, c+1, \ldots, c+n-1$.  We know that
the next to last tile $t_{d,j}(\alpha)$ contains the leftmost
remaining squares of rows $d,d+1, \ldots, d+n-1$, and so forth.
Thus, given our ordering, it is easy to explicitly determine all
the tiles of $\alpha$ in linear time.  Note also that this
algorithm yields an explicit proof of the fact (mentioned
earlier) that every acyclic partial ordering of a graph has an
acyclic extension.

\section{General Tilings And Abelian Height-Functions}
\label{generalsection}
\subsection{Examples}
The constructions in this paper allow us in principle to work out
the maximally informative abelian height function for any
polyomino tiling problem. The following facts, for example, are
simple to verify:

\begin{enumerate}
\item Given any coloring of the squares of $R$, a {\it generalized
ribbon tiling} is a tiling in which each tile contains exactly
one square of each color. Suppose, for example, we color the
square at position $(x,y)$ one of $mn$ colors according to the
values of $x \pmod{m}$ and $y \pmod{n}$.  In this case, the
canonical abelian height space is $\mathbb Z^{mn-1}$. (When
$m=n=2$, this is the square/skew-tetromino tiling problem
mentioned in \cite{Pr}.)
\item Let $T$ be the set of all polyominos with two black and two white
squares in the usual chess board coloring. (These are squares,
skew-tetrominos, and length-four vertical and horizontal bars.)
The canonical abelian height function space is $\mathbb Z \times
\mathbb Z/2\mathbb Z$.
\item If $T$ contains all trominos, then the canonical abelian height space is trivial.
\item Let $T$ be the set of all rectangular $m \times 1$ tiles and
$1 \times n$ tiles.  Then the canonical abelian height space is
$\mathbb Z^{(m-1)(n-1)}$.  (This problem was studied in \cite{KK},
where connectedness under local moves was proved using a
non-abelian height function.)
\item If $T$ contains only horizontal dominos, then the canonical abelian height space
is infinite dimensional.
\end{enumerate}

\subsection{Defining height functions at all vertices}
\begin{figure}
\begin{center}
\leavevmode \epsfbox[40 40 320 140]{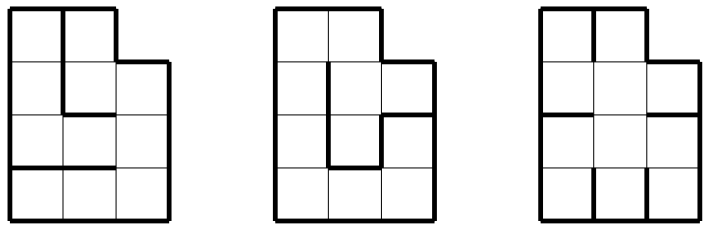}
\end{center}
\caption{Three possible spanning trees for a single tile.}
\label{spanningtrees}
\end{figure}

In general, the height function $h_{\alpha}$ is only defined at
the vertices of $R$ that are on the boundary of tiles of $\alpha$.
One (admittedly awkward) way around this is to replace each tile
with a spanning tree of the squares in that tile, and consider an
edge in $t$ to be on the boundary if the squares it divides are
not adjacent as members of the tree.  (See Figure
\ref{spanningtrees}.)  That is, if $v_1$ and $v_2$ are adjacent
vertices in $R$, we say $h_{\alpha}(v_2) - h_{\alpha}(v_1)$ is
equal to $\phi(v_1,v_2)$ when the squares incident to $(v_1,v_2)$
are either contained in different tiles or are contained in the
same tile $t$ but are not adjacent as members of the spanning
tree of $t$; this is enough to determine the value of
$h_{\alpha}$ at every vertex $v$ in $R$---and this is now true
even if we allow $T$ to contain non-simply-connected tiles.

Of course, the height values now depend on the particular choice
of spanning tree for each tile; one might even consider each way
of choosing a spanning tree of a tile $t$ to be a different tile
in its own right.  In this perspective, $T$ becomes a set of
tile--spanning-tree pairs.  Changing the spanning tree
corresponding to a given tile is a sort of local replacement move.

\subsection{Can distinct tilings have the same height function?}
We saw in the examples that the abelian height function space for
tromino-tilings is trivial; in this case, of course, all tilings
of $R$ by trominos have the same (trivial) height function
mapping all vertices to zero.

We will now describe a condition on tile sets $T$ that will ensure
that two distinct tilings do not have the same height function.
Given a map $\phi$ from the squares of $R$ to the abelian height
group $G$, we say a tile $t$ is {\it admissible} by $\phi$ if
$\phi(t) = 0$.  We say $t$ is {\it irreducible} with respect to
$\phi$ if there is no way to partition $t$ into two or more
smaller polyominos each of which is admissible. A tile set $T$ is
irreducible if each tile is irreducible with respect to the map
$\phi$ from the squares of $R$ to the canonical abelian height
space determined by those tiles.  Note that these definitions
depend only on the (canonically determined) values of $\phi$ on
squares $s$; they do not depend on any arbitrary choices we made
when choosing the value of $\phi$ on individual edges.

Now, consider two tilings $\alpha$ and $\beta$ of $R$ with the
same height function $h$.  We claim that each connected component
$u$ of the intersection of any tile in $\alpha$ with a tile in
$\beta$ satisfies $\phi(u) = 0$.  To see this, let $v_0, v_1,
\ldots, v_k=v_0$ be the vertices of a cycle surrounding $u$.
Since each edge in this cycle is a boundary edge of either
$\alpha$ or $\beta$, it follows that $$ \phi(u) = \sum_{s \in u}
\phi(s) = \sum_{i=1}^{k} \phi(v_i, v_{i-1}) = \sum_{i=1}^{k}
h(v_i) - h(v_{i-1}) = 0,$$ where the second sum is over clockwise
edges on the boundary of $u$. Unless $u$ is a tile in both
$\beta$ and $\alpha$, it must be contained in a larger tile $t$
in one of the two tilings, say $\beta$.  All the connected
components of intersections of $t$ with tiles of $\alpha$ are
admissible, and it follows that $t$ is reducible with respect to
$\phi$.  We can now deduce the following:

\begin{prop}
If $T$ is irreducible, then no two tilings of $R$ by tiles in $T$
have the same height function.
\end{prop}

This does not mean that height functions for reducible tile sets
cannot provide useful information, but if a tile set is reducible,
there may be aspects of the tiling that are not captured by the
height function.  A good example is the tile set $T$ containing
all domino tiles and all skew tetromino tiles. The abelian height
function for tilings by these tiles is simply the usual
domino-tiling height function, where a skew tetromino is treated
as a pair of dominos.

\subsection{Which functions are height functions of a tiling?}
Assume our height functions are defined at every vertex.  Then for
each oriented edge $e = (v_1, v_2)$ we define $\Delta_e$ to be the
set of all differences $h_{\alpha}(v_2) - h_{\alpha}(v_1)$ that
can possibly arise for a tiling $\alpha$ of any region $R$.
Recall that:
$$h_{\alpha}(v_2) - h_{\alpha}(v_1) = \phi(v_1,v_2) + \sum \phi(s),$$
where the latter sum ranges over all squares passed on the left
by the edge $(v_1,v_2)$.  Thus the set $\Delta_e$ is equal to
$\phi(v_1,v_2) + L_e$ where $L_e$ is the set containing $0$ and
all possible square sums $\sum \phi(s)$ over sets of squares that
can occur in the left side of a tile divided by $e$. Clearly,
$\Delta_{e} = - \Delta_{e^{-1}}$.

We say a function $h$ from the vertices of $R$ to $G$ has {\it
proper differences} with respect to the tile set if
\begin{enumerate}
\item $h(v_2) - h(v_1) = \phi(v_1,v_2)$ whenever $(v_1,v_2)$ is a
boundary edge of $R$.
\item $h(v_2) - h(v_1) \in \Delta_{(v_1,v_2)}$ whenever $(v_1,v_2)$ is an interior edge of $R$.
\item $h(v_0) = 0$ for some reference vertex $v_0$ on the boundary of $R$.
\end{enumerate}

We say a tile set $T$ is {\it complete} if it is irreducible and
for every simply connected region $R$, every function on $R$ with
proper differences is the height function of a tiling of $R$ by
tiles in $T$.  For this definition, $\phi$ is assumed to be the
usual map from edges into a group $G$ containing the canonical
abelian height space as a subgroup.  Like the definition of
irreducibility of $T$, this definition is independent of any
arbitrary choices we made in choosing $\phi$.  To see this, note
that if we add to $\phi$ a function $\phi'(e)$ on the edges of
$\mathbb Z^2$ such that $\phi'(s) = 0$ for any square, this
simply amounts to adding a tiling-independent function to each
$h_{\alpha}$ and adjusting the $\Delta_e$ values accordingly.
When a tile set is complete, there is a precise one-to-one
correspondence between tilings of $R$ by that tile set and height
functions on the vertices of $R$ with proper differences.

The fact that dominos are a complete set of tiles has been very
important to the study of domino tilings; it is often easier to
work with the space of height functions satisfying difference
restrictions than to work directly with domino tilings.
Completeness of domino tilings implies that existence algorithms
and random sampling algorithms for domino tilings can all be
derived from corresponding algorithms for spaces of height
functions with proper differences.

\subsection{Generalized Ribbon Tiling Completeness}
We show here that any set of generalized ribbon tilings is
complete.  Given a coloring of the squares of $\mathbb Z^2$ with
$n$ colors, recall that a generalized ribbon tile set $T$ consists
of all polyominos (together with corresponding spanning trees)
that contain exactly one square of each color.

Let $G$ be the canonical abelian height space for the generalized
ribbon tile set $T$.  Although we do not in general know the
structure of $G$, we can determine the structure of $H$, the
quotient of $G$ with the relations $\phi(s_1) - \phi(s_2) = 0$
whenever $s_1$ and $s_2$ have the same color.  Let $\psi$ be the
composition of $\phi$ with the quotient map sending $G$ to $H$
and observe that the value of $\psi(s)$ depends only on the color
of $s$ (the same is not necessarily true of $\phi(s)$). Let $e_i$
be the value of $\psi(s)$ when $s$ has color $i$.

Now $H$ is precisely the group generated by the $e_i$ subject to
the one relation $\sum_{i=0}^{n-1} e_i = 0$. (Because this
relation is enough to force $\psi(t)$ to be zero for each $t \in
T$, no additional relations are required.)  It is easy to see
that these tiles are irreducible: no nonempty proper subset of
the $e_i$ sums to zero in $H$, so each sum of $\phi(s)$ over a
nonempty proper subset of the squares of a tile has nonzero image
in $H$ and is hence a nonzero element of $G$.

Now, to show that $T$ is complete, we must show that every $h$
mapping the vertices of $R$ to $G$ with proper differences is the
height function of a tiling $\alpha$.  Suppose $h$ has proper
differences.  Given an edge $e=(v_1,v_2)$ in $R$, we define $b(e)
= h(v_2) - h(v_1) - \phi(e)$.  By the definition of proper
differences, $b(e)$ is the sum of $\phi(s)$ over some set $S(e)$
of squares $s$---and there exists at least one tile $t \in T$ such
that $e$ crosses $t$ and passes those squares on the left.  In
particular, this implies that $S(e)$ contains at most one square
of each color and cannot contain a square of every color.  Put
differently, the set $C(e)$ of colors that occur in squares of
$S(e)$ is a proper subset of $\{1,2,\ldots, n\}$, and the image
$\overline{b}(e)$ of $b(e)$ in $H$ is simply $\sum_{j \in C(e)}
e_j$.  We say $e$ has {\it color $j$ on its left} if $j \in
C(e)$.  Note that if $e$ has a square of color $j$ immediately to
its right and $h$ has proper differences, then $C(e)$ can never
contain $j$.

Now we will construct the tiling $\alpha$ for which $h$ is the
height function.  Let $a_1$, $a_2$, $a_3$, and $a_4$ be the
clockwise oriented edges of a square $s$ with color $i$. If $h$
is a function with proper differences and $v_0, v_1, v_2, v_3,
v_4=v_0$ are a path of vertices surrounding $s$ clockwise, we
have:
$$0 = \sum_{k=1}^{4} h(v_k) - h(v_{k-1}) = \sum_{k=1}^4 \phi(a_k) + b(a_k) = \phi(s) +
\sum_{k=1}^4 b(a_k)$$

Taking the image of the latter sum in $H$ gives
$$0 = e_i + \sum \overline{b}(a_k) = e_i + \sum_{k=1}^4 \sum_{j \in C(a_k)} e_j$$

Since none of the $C(a_k)$ can contain $i$; it is easy to see that
the only way $\sum \overline{b}(a_k)$ can be $-e_i$ in $H$ is if
the sum includes exactly one term of type $e_j$ for each $j \not
= i$.  Thus, for $1 \leq k \leq 4$, the $C(a_k)$ are disjoint, and
their union contains all colors except $i$.

Now, we say $s_1$ and $s_2$ are {\it tile adjacent} if they are
adjacent and the edge $e$ between them has a nonempty $S(e)$. We
define our tiles to be the connected components of the
tile-adjacency relation. If we can show that every tile contains
exactly one square of each color, we will have proved the
completeness of $T$.

Suppose $s_0$ has color $i$, and choose some $j \not = i$.  We
aim to show that $s_0$ is in the same tile as a square of color
$j$.  Observe that at least one edge of the four clockwise edges
of $S$ has color $j$ on its left. Let $s_1$ be the square adjacent
to $s_0$ that is incident to the edge $e$ for which this is the
case. Clearly, $s_1$ does not have color $i$.  If it has color
$j$, we are done.  Otherwise, suppose $a_1$, $a_2$, and $a_3$ are
the three other clockwise oriented edges of $s_1$.  We know that
the union of the $C(a_i)$ and $C(e^{-1})$ contains all colors
except the color of $s_1$.  Since $C(e^{-1}) \cup C(e) =
\{1,2,\ldots,n\}$, it follows that each of the $C(a_i)$ is a
proper subset of $C(e)$ and that one of them contains $j$. Now, we
let $s_2$ be the square adjacent to $s_1$ and incident to the
edge for which this is the case and repeat this process.  Since
the proper subsets cannot decrease in size indefinitely, we must
eventually reach a square of color $j$.  A similar argument shows
that such a path starting at $s_0$ never reaches another square
of color $i$.  It follows that each tile contains exactly one
square of each color and that the adjacency relation describes a
spanning tree of that tile.

Using similar arguments, one can show that the tile set containing
all $1 \times k$ and $m \times 1$ rectangles is complete.

\subsection{Completion of a Tile Set} A good example of a tile set
that is irreducible but not complete is, when $n\geq 3$, the set
of ribbon tiles of order $n$ with all but one of the possible
ribbon tile shapes.

When is it possible to add tiles to an incomplete tile set $T$ to
make it complete?  First, it is clear that we will have to add as
a tile every finite region $R$ on which there exists a height
function $h$ with proper differences such that the set $S$ of
edges $(v_1,v_2)$ in $R$ for which $h(v_2) - h(v_2) =
\phi(v_1,v_2)$ {\it does not} divide the squares of $R$ into more
than one component.  (There may be an infinite number of such
regions, and they may be arbitrarily large.)  Note that since all
of these new tiles are already admissible by $\phi$, adding these
tiles will not change the canonical abelian height space, and we
can use the same $\phi$ to define height functions for the new
tiling.

If each of these new tiles is irreducible, the result is a
complete set of tiles, called the {\it completion} of $T$.

\subsection{Non-Existence Proofs for Complete Tile Sets}
Conway-Lagarias \cite{CL} and Thurston \cite{T} observed that a
simple necessary condition for the existence of a tiling on $R$
is that the height function be well-defined on the boundary of
$R$. That is, the group-wise sum of $\phi(e)$ around the all the
edges $e$ of $R$ must be the identity.

It is not difficult to describe much stronger necessary
conditions.  If a tile set $T$ is complete and $R$ is a
simply-connected region, then (assuming the previous condition is
satisfied) the question of whether $R$ can be tiled by $T$
reduces to the question of whether there is a height function on
the vertices of $R$ with proper differences.  Even if $T$ is not
complete, the existence such a function is a necessary condition
for the existence of a tiling of $R$.  Thus, the following
interpolation problem is of interest:

Given a simply connected region $R$, a function $h$ defined on all
boundary vertices of $R$, and a set $\Delta_e$ for each oriented
edge (with $\Delta_e = -\Delta_{e^{-1}}$ for all $e$), is it
possible to interpolate $h$ to all vertices of $R$ in such a way
that $h(v_2) - h(v_1) \in \Delta_{(v_1,v_2)}$ whenever
$(v_1,v_2)$ is an edge in $R$?

When the height function space is a rank $n$ free-abelian group,
this is a special case of integer programming, which is in
general NP-complete.  But suppose we embed the height space in
$\mathbb R^n$ and let $\overline{\Delta_{e}}$ be the convex hull
of the set of points in $\Delta_e$.  Then a relaxed version of
the problem asks for an interpolation of $h$ to the interior
squares of $R$ where $h$ is allowed to assume any values in
$\mathbb R^n$ and the differences are merely required to lie in
$\overline{\Delta_e}$. Since each $\overline{\Delta_e}$ is a
convex polyhedron, this is now a linear programming problem that
can be solved in polynomial time.  If we prove that there is no
solution to the relaxed interpolation problem, we know $R$ cannot
be tiled by $T$.

\subsection{Note on Multiply Connected $R$} Suppose that $R$ is simply connected
except for $k$ ``islands,'' labeled $I_1, \ldots I_k$; we assume
that the boundary edge/vertex sets of the boundary of each $I_k$
and the outer boundary of $R$ are disjoint, but that each $I_i$ is
connected.

The problem with defining height functions on $R$ is that as we
traverse the boundary of a given $I_k$ clockwise and come back to
the same place, the height function changes by $\phi(I_i)$; thus,
if $\phi(I_i)\not = 0$, then the height function $h_{\alpha}$ is
not well-defined as a single-valued function. There are two ways
around this problem.  One is simply to accept that $h_{\alpha}$
is a multiple-valued function; only the differential $\nabla
h_{\alpha}(v_1,v_2) = h_{\alpha}(v_2) - h_{\alpha}(v_1)$, is well
defined, and this can be integrated to produce a single-valued
$h_{\alpha}$ on the vertices of any simply connected subset of
$R$---but not all of $R$.

Second, we could choose some replacement $\phi'$ for $\phi$
satisfying:
\begin{enumerate}
\item $\phi'(I_i)=0$ for each island $I_i$.
\item $\phi'(s) = \phi(s)$ for each $s$ in $R$.
\end{enumerate}

It is not hard to see that such a $\phi'$ always exists.  If we
define $h_{\alpha}$ using $\phi'$ instead of $\phi$, then it is
everywhere defined.  We require $h_{\alpha}(v_0)=0$ for some outer
boundary vertex $v_0$, and we also choose one vertex $v_i$ on each
island $i$. Clearly, the value of $h_{\alpha}$ on all boundary
vertices is determined by its value on the $v_i$.

In the special case of ribbon tilings, we define a {\it simple
local replacement move} on $R$ to be a move that involves removing
two tiles $t_1$ and $t_2$ that together comprise a {\it simply
connected} two-tile region $R'$ and replacing them with another
pair of tiles.  When $R$ is simply connected, no tiling of $R$
ever contains a pair $t_1$ and $t_2$ of tiles whose union is
connected but not simply connected---otherwise, there would be a
hole between them that could not be filled by order-$n$ ribbon
tiles.  Thus, when $R$ is simply connected, every possible local
replacement move is simple. It is not hard to show in general that
two order-$n$ ribbon tilings $\alpha$ and $\beta$ on $R$ can be
connected by a sequence of simple local replacement moves if and
only if $h_{\alpha}(v_i) = h_{\beta}(v_i)$ for each $i$. The
reader may also check that when $n<4$, every local replacement
move is simple.

\section{Open Problems}
\label{opensection} The domino tiling problem is perhaps the most
well understood of all tiling problems.  Height functions and
local replacement moves have played a fundamental role in the
proofs of many domino tiling theorems, including the following.

\begin{enumerate}
\item Proofs of local-move connectedness.
\item Algorithms for computation of minimal-length local move
sequences connecting two tilings.
\item Linear algorithms for determining when a region $R$ has a tiling.
\item Mixing time bounds on random processes based on local moves. \cite{LRS}, \cite{RT}
\item Perfect sampling algorithms based on ``coupling from the past.'' \cite{PW}
\item Asymptotic entropy computations for random tilings and descriptions
of finite and asymptotic Gibbs measures. \cite{CEP}, \cite{CKP}
\item Law of large numbers and large deviations principles for
asymptotic height function shapes of random tilings. \cite{CKP}
\item Conformal invariance of certain tiling properties in the asymptotic
scaling limit. \cite{K2}
\end{enumerate}

These results and references are ridiculously far from
exhaustive; there have been numerous papers relating to each of
the above topics and many others.  We include the references not
to assign credit, necessarily, but to provide pointers to papers
with more complete bibliographies.  Because ribbon tilings are in
many ways the simplest polyomino tilings that generalize domino
tilings, the following question is very natural:

\begin{question} Which of the above results can be generalized to ribbon
tilings?
\end{question}

An additional task is to generalize to other tiling spaces with
abelian height functions.

\begin{question}
Can connectedness under local replacement moves be generalized to
the other tiling problems in this paper?  In particular, is there
a local move result for every complete, translation-invariant
tile set?
\end{question}

Finally, the author hopes this paper will be a step towards a
more thorough understanding of at least the abelian case of the
Conway-Lagarias height function theory.  The following is a more
theoretical question:

\begin{question}
Is there a simple classification of the complete, translation
invariant tile sets with finitely many tile shapes?  Given such a
tile set $T$, is there always a polynomial algorithm for
determining whether a simply connected region can be tiled by $T$?
\end{question}

\section{Acknowledgements}
Many thanks to Henry Cohn for bringing the problem to my
attention, for numerous helpful conversations, and for reviewing
early drafts of the paper.  Thanks also to L\'{a}szl\'{o}
Lov\'{a}sz for helpful conversations and to Igor Pak for
suggesting additional references.


\end{document}